\documentclass{amsart}
\usepackage[dvips]{graphicx}
\usepackage{amsmath,graphics}
\usepackage{amsfonts,amssymb}
\usepackage{comment}
\usepackage{yfonts}

\theoremstyle{plain}

\newtheorem{theorem}{Theorem}[section]
\newtheorem{lemma}[theorem]{Lemma}
\newtheorem{corollary}[theorem]{Corollary}
\newtheorem{proposition}[theorem]{Proposition}

\newtheorem{remark}[theorem]{Remark}
\newtheorem{remarks}[theorem]{Remarks}

\newtheorem*{notations}{Notations}
\newtheorem{example*}{Example}

\def\Q{\mathbb{Q}} 
 
\def\Z{\mathbb{Z}} 
\def\R{\mathbb{R}} \def\C{\mathbb{C}}
\def\N{\mathbb{N}} \def\P{\mathbb{P}}

\begin{document}
\title[Fibered cohomology classes in dimension three...]{Fibered cohomology classes in dimension three, twisted Alexander polynomials and Novikov homology}
\author{Jean-Claude Sikorav}
\address{Unit\'e de Math\'ematiques Pures et Appliqu\'ees, UMR CNRS 5669, \'Ecole normale sup\'erieure de Lyon\\   France}
\email{jean-claude.sikorav@ens-lyon.fr}

\date{May 10, 2021}


\begin{abstract}
We prove that for ``most'' closed $3$-dimensional manifolds $M$, the existence of a closed non singular one-form in a given cohomology class $u\in H^1(M,\R)={\rm Hom}(\pi_1(M),\R)$ is equivalent to the fact that every twisted Alexander polynomial $\Delta^H(M,u)\in\Z[G/\ker u]$ associated to a normal subgroup with finite index $H<\pi_1(M)$ has a unitary $u$-minimal term. 
\end{abstract}

\maketitle

\section{Introduction and statement of the main result}
We consider $M$ a closed connected $3$-manifold. Let $G:=\pi_1(M)$ and let
$u$ be a nonzero element of ${\rm Hom}(G,\R)$, which will be identified with $H^1(M,\R)$.
Denote by ${\rm rk}(u)$ the rank of $u$, i.e. the number of free generators of $G/\ker u$. We are interested in the following

\smallskip
\noindent{\bf Question.} {\it Does there exist a nonsingular closed $1$-form $\omega$ in the class $u$?}

\smallskip
If such a form exists, we say that $u$ is {\it fibered}. The reason is that if ${\rm rk}(u)=1$ so that $au(G)\subset\Z$ for a suitable $a\ne0$, such a form is $a^{-1}f^*dt$ for $f$ a fibration to $S^1=\R/\Z$. More generally, by \cite{Tischler 1970}, if $u$ fibers then $M$ fibers over $S^1$: perturb $\omega$ to $\omega'=\omega+\varepsilon$ such that ${\rm rk}([\omega'])=1$ and $\varepsilon$ is $C^0$-small. Then $\omega'$ is still nonsingular, thus $M$ fibers.

\smallskip
An answer to this question was given in rank one by \cite{Stallings 1962}: {\it if ${\rm rk}(u)=1$, $u$ fibers if and only if $\ker u$ is finitely generated}. Actually, Stallings required $M$ to be irreducible, but using Perelman it is unnecessary. 

\smallskip
In any rank, the paper \cite{Thurston 1986} introducing the Thurston (semi-)norm on $H^1(M;\R)$ proved the following results: 1) the unit ball of the norm is an ``integer polyhedron'', i.e. it is defined by a finite number of inequalities $u(g)\le n$, $g\in G$, $n\in\N^*$; 2) {\it the set of fibered $u\in H^1(M;\R)\setminus\{0\}$ is a cone over the union of some maximal open faces of the unit sphere of the Thurston norm}. Note that thanks to Stallings, to know if a given face is ``fibered'', it suffices to test one element $u$ of rank one and see if $\ker u$ is finitely generated.

\smallskip
In the 2000s and beginning of 2010s, S. Friedl and S. Vidussi studied this question again, mostly in rank one, in connection with what was then a conjecture of Taubes: {\it $u$ fibers if and only if $u\wedge [dt]+a\in H^1(M\times S^1)$ is represented by a symplectic form, where $a\in H^2(M;\R)$ satisfies $a\wedge u\ne0$}. The starting point was the relation of Seiberg-Witten invariants of $M\times S^1$ and {\it twisted Alexander polynomials}, see below and Section \ref{twistedalexanderpolynomials}.
They ultimately solved that conjecture in \cite{Friedl-Vidussi 2013}, and obtained as a byproduct a new answer for the characterization of fibered classes in the case of rank $1$: {\it if ${\rm rk}(u)=1$, $u$ fibers if and only if 
all twisted Alexander polynomials $\Delta^H(G,u)$ are nonzero.}

\smallskip
Let us describe briefly what are these twisted Alexander polynomials (for a detailed presentation, see \cite{Friedl-Vidussi 2011}). In fact, we do it only for a special case, which is already sufficient: those associated to finite covers, see \cite{Friedl-Vidussi 2008}, section 3.2. 

\smallskip
Recall first the definition of the {\it order} of a finitely generated module ${\mathcal M}$ over a Noetherian UFD $R$: it is the greatest common divisor of the $p$-minors of $A$ in a finite presentation
$$R^q\stackrel{\times A}{\longrightarrow} R^p\to {\mathcal M},$$ 
where $\times A$ is the right multiplication by a  matrix $A\in{\rm M}_{q,p}(R)$. Thus it is an element of $R$ defined up to multiplication by a unit. We denote it by ${\rm ord}_R({\mathcal M})$, usually viewed as an element of $R$. See Section 3.1.

\smallskip
 Since $G/\ker u\approx\Z^r$, the ring $\Z[G/\ker u]$ is isomorphic to 
 $\Z[t_1^{\pm1},\cdots,t_r^{\pm 1}]$,
 thus it is a Noetherian UFD. In particular, this order vanishes if and only if there are no $p$-minors or they all vanish.

 \smallskip
Then let $H$ be a normal subgroup of $G$ with finite index, denoted by $H\triangleleft_{f.i.} G$. We define $H_1(H;\Z[G/\ker u])$ as the homology of $H$ with coefficients in the $H$-module $\Z[G/\ker u]$. It is naturally a module over $\Z[G/\ker u]$ by action on the coefficients, which is finitely generated since $H$ is finitely generated. By definition, the twisted Alexander polynomial of $(G,u)$ associated to $H$ is
 $$\Delta^H(G,u):={\rm ord}_{\Z[G]}(H_1(H;\Z[G/\ker u])).$$
 Since the units of $\Z[G/\ker u]$ are $\pm (G/\ker u)$ (i.e. $\pm t_1^{i_1}\cdots t_r^{i_r})$, it is an element of $\Z[G/\ker u]/{\pm (G/\ker u)}$.

\smallskip
 It is not too difficult to prove that, if $u$ is fibered, $\Delta^H(G,u)$ is always {\it $u$-monic}, i.e. its $u$-minimal term has a coefficient $\pm1$: see Proposition \ref{fiberingimpliesalexanderismonic}. In the rank one case, this goes back to Alexander.

\medskip
We can now state our main result.

\begin{theorem} \label{mainresult} Let $M$ be a closed $3$-manifold such that $\widetilde M$ is contractible and $G:= \pi_1 (M)$ is virtually residually torsion-free nilpotent (VRTFN), and let $u$ be a nonzero element of ${\rm Hom}(G,\R)=H^1(M,\R)$. 

\smallskip
Assume that  $\Delta^H (G,u)$ is $u$-monic for every  $H\triangleleft_{f.i.} G$. Then $u$ is fibered, i.e. represented by a nonsingular closed $1$-form.
  \end{theorem}
  
\subsection*{Comments.}   1) Building on \cite{Agol 2014}, \cite{Koberda 2013} proves that $\pi_1(M)$ is VRTFN for all geometric manifolds which are not Sol.  In particular, if $M$ is hyperbolic, this follows from the fact that $\pi_1(M)$ is virtually a right-angled Artin group.

\smallskip
If $M$ is Sol, $\pi_1(M)$ is not virtually nilpotent, but $M$ is either a torus bundle over $S^1$ with hyperbolic monodromy or 
has a finite cover of this type and $H^1(M,\R)=0$. Thus in that case the theorem is obvious.

\smallskip
The hypothesis that $\widetilde M$ is contractible can be dispensed with: if it does not hold, then (since $b_1(M)>0$) we are in one of the two following cases: either $M$ is nonprime thus nonfibered and the twisted Alexander polynomials always vanish; or $M$ fibers over $S^1$ with fiber $S^2$ or $\R\P^2$.

\medskip
2) In rank one, our result is weaker than \cite{Friedl-Vidussi 2013}. However, even in that case we believe that our proof, which is based on different ideas, may be of interest.

\section{Sketch of the proof and content of the paper}\label{sketchoftheproof}

The main idea is to express the fibering condition on $u$ by the vanishing of some {\it Novikov homology} associated to $(G,u)$ and the nonvanishing of $\Delta^H(G,u)$ by the vanishing of some {\it Abelianized relative Novikov homology} associated to $(G,H,u)$. 

\smallskip
In turn, these vanishings are expressed by the invertibility of some matrix in the Novikov ring associated to $(G,u)$ and of its image in the Novikov ring associated to $(G/(H\cap\ker u),\overline u)$ ($\overline u$ being induced by $u$).

\smallskip
Then the theorem is reduced to a result about ``finite detectability of invertible matrices'' for a VRTFN group.

\medskip
We now describe the content of the paper.

\medskip
In Section \ref{twistedalexanderpolynomials}, we define the twisted Alexander polynomials $\Delta^H(G,u)$.

\smallskip
In Section \ref{novikovhomology}, we define the Novikov ring $\Z[G]_u$, and the Novikov homology $H_*(G,u)$. We quote the result of \cite{Bieri-Neumann-Strebel 1987}, building upon previous results of Stallings and Thurston: {\it if $G=\pi_1(M)$ with $M$ a closed $3$-manifold, $u$ fibers if and only if $H_1(G,u)=0$}.

\smallskip
In Section \ref{abelianizednovikov}, we explain the relations between twisted Alexander polynomials and an ``Abelian relative'' version of Novikov homology.

\smallskip
In Section \ref{computationsindimensionthree}, we specify the computation of $H_1(G,u)$ for $G = \pi_1(M^3)$ with $\widetilde M$ contractible, thanks to the form of a presentation of $G$ given by a Heegaard decomposition and Poincar\'e duality. We deduce that ($H_1(G,u)=0)$ is equivalent to the invertibility in ${\rm M}_{p-1}(\Z[G]_u)$ of some matrix $A  \in {\rm M}_{p-1}(\Z[G])$ where $p$ is the genus of the decomposition. 

\smallskip
Similarly, the vanishing of $H_1^{ab}(G/(H\cap\ker u)),\overline u)$ is equivalent to the invertibility of the image of $A$ in ${\rm M}_{p-1}(\Z[G/H \cap  \ker u]_{\overline u})$. Thus we have reduced Theorem \ref{mainresult} to  Theorem \ref{mixedfinitedetectability}: 

\smallskip
{\it If $G$ is finitely generated and VRTFN, a matrix $A  \in {\rm M}_n(\Z[G])$ whose image in ${\rm M}_n(\Z[G/H]_u)$ is invertible for every $H\triangleleft_{f.i.} G$, is invertible in ${\rm M}_n(\Z[G]_u)$}. 

\smallskip
\noindent The restriction to finitely generated groups is actually not necessary, but simplifies the proof.

\medskip
Theorem \ref{mixedfinitedetectability} is proven in Section \ref{proofofmixedfinitedetectability}. There are two main ingredients:  

\begin{itemize}

\item (Sections \ref{finitelydetectableunitsandfullleftideals} to \ref{section:theorem7.3fornilpotent}) the case when $G$ is nilpotent, which uses three key facts: 1)  \cite{Hall 1959} a simple $\Z[G]$-module is finite;  2) \cite{Goldie 1958} when $G$ is nilpotent and torsion-free, $\Z[G]$ has a classical ring of quotients on the right; 3) \cite{Chatters 1984} when $G$ is nilpotent and torsion-free, $\Z[G]$ is a UFD in the non commutative sense.

\smallskip
\item (Section \ref{section:malcevneumanncompletion}) the fact that when $G$ is RTFN, it is orderable, thus one can embed $\Z[G]$  in the Mal'cev-Neumann completion $\Q\langle G \rangle$, which is a division ring (or skew field); moreover, by a remark of \cite{Kielak 2020} the order can be chosen so that $\Q\langle G \rangle $ contains $\Z[G]_u$. Actually, we work mostly with a subfield introduced by \cite{Eizenbud-Lichtman 1987}, which contains $\Z[G]$ and whose elements have ``controlled'' support.

\end{itemize}

\begin{notations} {\rm In the following text, $G$ is a finitely generated group and $u:G\to\R$ a nonzero homomorphism. Thus 
$G/\ker u \approx \Z^r$, $r={\rm rk}(u)$, 
and 
$\Z[G/\ker u] \approx  \Z[t_1^{\pm},\cdots,t_r^{\pm1}]$ is a UFD (unique factorization domain). }
\end{notations}

\section{Twisted Alexander polynomials}
\label{twistedalexanderpolynomials}

\medskip
\subsection{Order of a finitely generated module over a Noetherian UFD}  Let ${\mathcal M}$ be a finitely generated $R$-module where $R$ is a Noetherian UFD. One defines (cf. \cite{Eisenbud})

\begin{itemize}

\smallskip
\item the  {\it Fitting ideal (or elementary)  of order $0$} ${\rm Fitt}_0({\mathcal M})$ as the ideal of $R$ generated by  the  $p$-minors of a matrix $A\in{\rm M}_{q,p}(R)$ where 
$$R^q\stackrel{\times A}{\longrightarrow} R^p\stackrel{p}{\longrightarrow}  {\mathcal M}$$
is a presentation of ${\mathcal M}$ with $\times A$ the multiplication on the right by $A$;

\smallskip
\item
the {\it order} ${\rm ord}_R({\mathcal M})$ as the greatest common divisor (gcd) of ${\rm Fitt}_0({\mathcal M})$. 

\end{itemize}

\begin{remark} 
{\rm We use multiplication on the right rather than on the left since later we will have mostly noncommutative rings, and we prefer to work with left modules. Thus elements of $R^q$, $R^p$ are interpreted as row vectors.}

\end{remark}

\smallskip
It is easy to prove that the definition of  ${\rm Fitt}_0({\mathcal M})$ and thus of ${\rm ord}_R({\mathcal M})$  does not depend on the presentation: if $R^{q_1}\stackrel{\times A_1}{\longrightarrow} R^{p_1}\stackrel{\pi_1}{\longrightarrow}{\mathcal M}$ and $R^{q_2}\stackrel{\times A_1}{\longrightarrow} R^{p_2}\stackrel{\pi_2}{\longrightarrow}{\mathcal M}$ are two presentations, one can lift $\pi_1$ to $\times B:R^{m_1}\to R^{m_2}$ and obtain a presentation 

$$R^{q_1}\oplus R^{q_2}\stackrel{\times A}{\longrightarrow} R^{p_1}\oplus R^{p_2}\stackrel{\pi_1+\pi_2}{\longrightarrow} {\mathcal M}\ , \ A=\begin{pmatrix}{\rm I}_{p_1}&B\\ 
0&A_2
\end{pmatrix}.$$
Similarly, there is a presentation

$$R^{q_1}\oplus R^{q_2}\stackrel{\times C}{\longrightarrow} R^{p_1}\oplus R^{p_2}\stackrel{\pi_1+\pi_2}{\longrightarrow}{\mathcal M}\ , \ C=\begin{pmatrix}A_1&0\\ 
D&{\rm I}_{p_2}
\end{pmatrix}.$$
Thus the kernel of $\pi_1+\pi_2$ is the row space of $\begin{pmatrix}{\rm I}_{p_1}&B\\ 
0&A_2
\end{pmatrix}$, and also that of $\begin{pmatrix}A_1&0\\ 
D&{\rm I}_{p_2}
\end{pmatrix}$: this implies (for any ring) that the ideals generated by the $p_1$-minors of $A_1$ and by the $p_2$-minors of $A_2$ coincide. 

\smallskip
The main property of this order is the

\begin{proposition}
\label{orderandannihilator} Let $R$ be a Noetherian UFD, ${\mathcal M}$ an $R$-module generated by $p$ elements, and ${\rm ann}_R({\mathcal M})$ its annihilator. Then one has the divisions
$${\rm gcd}({\rm ann}_R({\mathcal M})) \mid {\rm ord}_R({\mathcal M}) \mid \big({\rm gcd}({\rm ann}_R({\mathcal M}))\big)^p.$$
More precisely, one has the inclusions of ideals
$${\rm ann}_R({\mathcal M})\supset{\rm Fitt}_0({\mathcal M})\supset({\rm ann}_R({\mathcal M}))^p=\langle a^p\mid a\in{\rm ann}_R({\mathcal M})\rangle$$
\end{proposition}

\noindent{\it Proof.} Let $R^q\stackrel{\times A}{\longrightarrow} R^p\stackrel{p}{\longrightarrow}{\mathcal M}$ be a presentation of ${\mathcal M}$, with $A\in{\rm M}_{q,p}(R)$. Then 
$$a\in{\rm ann}_R({\mathcal M})\Leftrightarrow aR^p\subset R^qA\Leftrightarrow(\exists X\in{\rm M}_{p,q}(R) )\  XA=a{\rm I}_p.$$

Let $\mu$ be a $p$-minor of $A$. Changing the order of the coordinates, we have $A=\begin{pmatrix}A_1\\A_2\end{pmatrix}$ where $A_1\in {\rm M}_p(R)$ with $\det A_1=\mu$. Thus there exists $X_1\in {\rm M}_p(R)$ such that $X_1A_1=\mu{\rm I}_p$, and $X=\begin{pmatrix}X_1&0\end{pmatrix}\in{\rm M}_{p,q}(R)$ satisfies $XA=\mu{\rm I}_p$. Thus ${\rm ann}_R({\mathcal M})$ contains $\mu$, thus it contains ${\rm Fitt}_0({\mathcal M})$.

\medskip
For the right inclusion: by the Cauchy-Binet formula for $\det(XA)$, the identity $XA=a{\rm I}_n$ implies that $a^p=\sum \mu_{X,i}\mu_{A,i}$ where $\mu_{X,i}$ and $\mu_{A,i}$ are $p$-minors of $A$ and $X$, thus $a^p\in{\rm Fitt}_0({\mathcal M})$. $\qed$

\subsection{Twisted Alexander polynomials.}
\label{twistedalexander} For every subgroup $H\triangleleft_{f.i.}G$, consider the left $G$-module 
$\Z[G/\ker u]^{G/H},$
where $G$ acts naturally both on $G/\ker u$ and on $G/H$, thus permuting the factors $\Z[G/\ker u]$. Thus one can define the homology group
$$H_1(G;\Z[G/\ker u]^{G/H}),$$ 
which is a module over $\Z[G]$ or over $\Z[G/(H\cap\ker u)]]$, but not on $\Z[G/\ker u]\approx\Z[G/\ker u]$. 

\smallskip
On the other hand, the action of $G$ on $G/\ker u$ descends to an action of $G/\ker u$ on $\Z[G/\ker u]^{G/H}$ which does not permute the factors $\Z[G/\ker u]$, and $H_1(G;\Z[G/\ker u]^{G/H})$ becomes a module over $\Z[G/\ker u]$. 

\smallskip
Viewing $H_1(G;\Z[G/\ker u]^{G/H})$ as a module over $\Z[G/\ker u]$, which is isomorphic to $\Z[\Z^r]$ thus still a UFD, we define
the {\it twisted Alexander polynomial associated to $H$}:
$$\Delta^H(G,u):={\rm ord}_{\Z[G/\ker u]}(H_1(G/\ker u;\Z[G]^{G/H})\in \Z[G/\ker u]\ {\rm mod}\ \pm G.$$

\medskip
\begin{remarks}

 {\rm 
\begin{enumerate}

\smallskip
\item
When $H=G$, $\Delta^G(G,u)$ can be denoted by $\Delta(G,u)$ and called ``multivariate Alexander polynomial'',  related to the Alexander polynomial of links. If ${\rm rk}(u)=1$, one recovers the classical Alexander polynomial, as generalized by \cite{Milnor 1968}.

\smallskip
\item If $X$ is a finite complex with $\pi_1(X)=G$ and $\widehat X_{H,u}$ the covering associated to $\ker u$, $H$ and $H\cap\ker u$, we have an isomorphism of modules over $\Z[H/H\cap\ker u]\approx\Z[u(H)]$:
$$H_1(G;\Z[G/\ker u]^{G/H})\approx H_1(\widehat X_{u,H};\Z).$$
One can deduce that 
$$\Delta^H(G,u)=0\ {\rm is}\ u-{\rm monic}\ \Leftrightarrow\Delta(H,u_{|H})\ {\rm is}\ u-{\rm monic}.$$

\end{enumerate}
}
\end{remarks}

\subsection{Comparison with \cite{Friedl-Vidussi 2008}{\rm, 3.2.1 to 3.2.4}.} (I change their notation from $N$ to $M$).
They start from 

\smallskip
$\bullet$ a free Abelian group $F$ together with a morphism $\psi:G=\pi_1(M)\to F$: in our case, $F=G/\ker u$ and $\psi$ is the natural projection. 

\smallskip
$\bullet$ a morphism $\gamma:G\to{\rm GL}(k,\Z[F])$: in our case, this is the morphism $G\mapsto{\rm GL}(\Z[G/\ker u])^{G/H}$ induced by the action of $G$ on $G/H$ but not on $G/\ker u$.

\smallskip
\noindent Thus their $\alpha=\gamma\otimes\psi$ is the morphism $G\mapsto{\rm GL}(\Z[G/\ker u])^{G/H}$ induced by the actions of $G$ on $G/\ker u$ and on $G/H$.

\smallskip
Then they define for any $\alpha:G\to{\rm GL}(k,\Z[F])$ the {\it i-twisted Alexander polynomial of $(G,\alpha)$}, denoted by $\Delta^\alpha_{M,i}$, by 
$$\Delta^\alpha_{G,i}={\rm ord}_R(H_i(M;\Z[F]^k))$$ where the (hidden) action of $G$ on $\Z[F]^k$ is $\alpha$. In our notations, we thus have
$$H_i(G;\Z[G/\ker u]^{G/H})=H_i(M;\Z[F]^k).$$
(for all $i$ if $\widetilde M$ is contractible, for $i\le 1$ if not).
Thus
$$\Delta^H(G,u)=\Delta^{\gamma\otimes\psi}_{M,1}.$$

\subsection{Fibering implies $u$-monicity of Alexander polynomials}

\begin{proposition}

\label{fiberingimpliesalexanderismonic}
If $G=\pi_1(M)$, $M$ a closed manifold of any dimension and $u\in H^1(M,\R)$ fibers, then $\Delta^H(G,u)$ is always $u$-monic.

\end{proposition}

\noindent{\it Proof.} Let $\omega$ be a nonsingular form in the class $u$, and let $X$ be a vector field on $M$ such that $\omega(X)=1$. On the universal cover $\widetilde M$, $\omega$ lifts to $df$, $X$ lifts to $\widetilde X$ with $df(\widetilde X)=1$. Thus the flow $(\varphi_X^t)$ lifts to $(\varphi_{\widetilde X}^t)$  with $f\circ\varphi_{\widetilde X}^t-f=t$. 

\smallskip
Fix a small cell decomposition of $M$ and lift it to $\widetilde M$, and choose lifts $\sigma\mapsto\widetilde\sigma$ of cells in $M$. For $t>0$ large enough, $\varphi_{\widetilde X}^t$ is equivariantly homotopic to an equivariant chain map $\widetilde\varphi$ such that $\widetilde\varphi(\widetilde\sigma)=g(\widetilde\tau)$ with $u(g)>0$. 

\smallskip
Thus the identity of the cell complex $C_*(\widetilde M)$ is homotopic over $\Z[G]$ to a chain map $A$ whose support in $G$ lies in $\{u>0\}$. 
Thus $A$ induces the identity on $H_1(M,\Z[u(G)]^{G/H})$. View $A$ as a matrix in some ${\rm M}_N(\Z[G])$, and denote by $\overline A$ its image in ${\rm M}_{N[G:H]}(\Z[u(G)])$ which acts on $C_*(\widetilde M)\otimes_{\Z[G]}\Z[u(G)]^{G/H})$. 

Then the support of $\overline A$ in $u(G)$ lies in $]0,+\infty[$. On the other hand, since $A$ induces the identity on $H_1(M,\Z[u(G)]^{G/H})$
 $\det({\rm Id}-\overline A)$ annihilates $H_*(M,\Z[u(G)]^{G/H})$ and in particular $H_1(M,\Z[u(G)]^{G/H})$. Since $\det({\rm Id}-A)$ is $u$-monic, we are done.

\section{Novikov homology} \label{novikovhomology}

\subsection{Novikov ring} We define the {\it Novikov ring} $\Z[G]_u$ as the following group of formal series over $G$ with coefficients in $\Z$:
$$\Z[G]_u :=\{\lambda \in \Z[[G]]\mid ( \forall C\in\R)\  {\rm supp}(\lambda)\cap \{u\le C\}\ \hbox{\rm is finite}\},$$
where $\{u\le C\}=\{g\in G\mid u(g)\le C\}$. It is easy to see that the multiplication
$$\sum_{g_1\in G} a_{g_1}g_1.\sum_{g_2\in G} b_{g_2}g_2= \sum_{g \in G} \big(\sum_{g_1g_2=g} a_{g_1}b_{g_2}\big)g$$
is well defined and makes $\Z[G]_u$ a ring containing $\Z[G]$ as a subring.

\medskip
\noindent{\bf Units of $\Z[G]_u$.} If $\lambda=1+a\in \Z[G]_u$ and ${\rm supp}(a)\subset\{u>0\}$, $\lambda$ is invertible, with inverse $\lambda^{-1}=\displaystyle\sum_{n=0}^\infty(-a)^n$. Thus every element of $\Z[G]_u$ whose $u$-minimal part is of the form $\pm g$ with $g\in G$, is a unit. We call such elements {\it $u$-{monic}}. More generally, if $A\in{\rm M}_n(\Z[G]_u)$ with ${\rm supp}(A)\subset\{u>0\}$, ${\rm I}_n+A$ is invertible with $({\rm I}_n+A)^{-1}=\displaystyle\sum_{n=0}^\infty(-A)^n$.

\smallskip
In the case when $\Z[G]$ has no zero divisors, in particular for $G=\Z^r$, every unit of $\Z[G]_u$ is $u$-monic, thus the units of $\Z[G]_u$ coincide with $u$-monic elements.

\medskip
\subsection{Novikov homology, relation with fibering} The Novikov homology $H_*(G,u)$  is defined as the homology of $G$ with coefficients in the left $\Z[G]$-module $\Z[G]_u$:
$$H_*(G,u):= H_*(G, \Z[G]_u).$$ 
Although we shall not use it explicitly, let us quote the following easy result, which was a great inspiration for our work. Note the relation with \cite{Stallings 1962}.

\begin{theorem} \cite{Bieri-Neumann-Strebel 1987}, \cite{Sikorav 1987}
\label{finitelygeneratedkernel} If ${\rm rk}(u)=1$, the kernel of $u$ is finitely generated if and only if $$H_1(G,u)=0=H_1(G,-u).$$
\end{theorem}

\smallskip
For this paper, the interest of Novikov homology lies in the following

\begin{theorem} {\bf \cite{Bieri-Neumann-Strebel 1987}.} Let $G=\pi_1(M)$, where $M$ is a closed and connected three-manifold. The following are equivalent: 
\begin{itemize}
\item $u$ is fibered.

\smallskip
\item $H_1(G,u) = 0$.

\end{itemize}
\end{theorem}

\begin{remarks}  {\rm \

\smallskip
\begin{enumerate}

\smallskip
\item This is their Theorem E, reinterpreted in terms of Novikov homology, cf. p.456 of the paper.

\smallskip
\item At the time, one needed $M$ to contain no fake cells, and also the hypothesis $\pi_1(M)\ne\Z\oplus\Z/2\Z$ (to avoid a possible fake $\R\P^2\times S^1$), restrictions removed later thanks to Perelman.

\smallskip
\item  In \cite{Sikorav 1987}, the equivalence between (i) and ($H_1(G, u)=0=H_1(G,-u)$) was proved as an immediate consequence of \cite{Stallings 1962}, \cite{Thurston 1986} and Theorem \ref{finitelygeneratedkernel}. This would suffice to prove our main result with almost no change.

\end{enumerate}

}

\end{remarks}

\subsection{Computation of $H_1(G, u)$}

To simplify the notations, we assume that $G$ is finitely presented (anyhow, we only need this case). Let $\langle x_1,\cdots,x_p\mid r_1,\cdots,r_q \rangle$ be a presentation of $G$, and let $D_1$, $D_2$ be defined as in Section \ref{twistedalexanderpolynomials}. Then, denoting by $(D_i)_u\in{\rm M}_{q,p}(\Z[G]_u)$ the matrix obtained by the base change $\Z[G] \to \Z[G]_u$, we have 
$$H_1(G, u) = \ker(\times (D_1)_u) /{\rm im}(\times D(_2)_u).$$
Since $u \ne 0$, there exists $i\in\{1,\cdots,p\}$, such that $u(x_i)\ne 0$, thus $x_i-1$ is invertible in $\Z[G]_u$. We can assume that $i=p$. Denoting by $D_2^{(i)}$ the matrix obtained by deleting the $i$-eth column of $D_2$, we have $H_1(G, u) \approx {\rm coker} (\times D_2^{(i)})_u$.

\begin{corollary} 
\label{deletingonecolumn}
We have
$$H (G,u) = 0 \Leftrightarrow \times (D_2^{(i)})_u : \Z[G]_u^q \to \Z[G]_u^{p-1}\ \hbox{\rm is onto}.$$
\end{corollary}
Equivalently, there exists $\widetilde X\in{\rm M}_{p-1,q}(\Z[G]_u)$ such that $\widetilde XD_2^{(i)}={\rm I}_{m-1}$. 
By truncating $\widetilde X=\displaystyle\sum_{g\in G}X_gg$, $X_g\in{\rm M}_{p-1,I}(\Z)$ below a sufficiently high level of $u$, i.e. by defining the finite sum $X=\displaystyle\sum_{u(g)\le C}X_gg$ with $C$ sufficiently large, we obtain $X\in{\rm M}_{p-1,I}(\Z[G])$ such that  $XD_2^{(i)}={\rm I}_{p-1}+A$ with $u>0$ on ${\rm supp}(A)$. Since such a matrix is invertible over $\Z[G]_u$, we obtain the  following

\begin{corollary} 
\label{h1=0withoutnovikovring}
($H_1(G,u;R) = 0$) is equivalent to the existence of a matrix $X\in{\rm M}_{m-1,I}(R[G])$  such that $XD_2^{(i)}={\rm I}_{m-1}+A$ with $u>0$ on ${\rm supp}(A)$.
\end{corollary}

\begin{corollary} 
\label{characterizationoverzg}
We have
$$H_1(G,u) = 0 \Leftrightarrow (\exists X\in{\rm M}_{p-1,q}(\Z[G]))\ u>0 \ {\rm on}\ {\rm supp}(XD_2^{(i)}-{\rm I}_{p-1}).$$
\end{corollary}

\noindent{\it Proof.} By Proposition \ref{deletingonecolumn}, the left hand side is equivalent to the existence of  $\widetilde X\in{\rm M}_{p-1,q}(\Z[G]_u)$ such that $\widetilde X(D_2^{(i)})={\rm I}_{p-1}$. 
By truncating $\widetilde X=\displaystyle\sum_{g\in G}X_gg$, $X_g\in{\rm M}_{p-1,q}(\Z[G])$ below a sufficiently high level of $u$, i.e. by defining the finite sum $X=\displaystyle\sum_{u(g)\le C}X_gg$ with $C$ sufficiently large, we obtain $X\in{\rm M}_{p-1,q}(\Z[G])$ such that $u>0$ on ${\rm supp}(XD_2^{(i)}-{\rm I}_{p-1})$. Conversely, since such a matrix $X$ is invertible over $\Z[G]_u$, this proves the corollary.

\section{Abelianized relative Novikov homology and twisted Alexander polynomials}
\label{abelianizednovikov}

Consider the induced morphism
$\overline u:G/\ker u\to\R$ and the associated Novikov ring $\Z[G/\ker u]_{\overline u}$, which is a left $\Z[G]$-module, and define the {\it  Abelianized Novikov homology}
$$H_1^{ab}(G,u):=H_1(G;\Z[G/\ker u]_{\overline u}).$$ 
It is in fact the original homology defined in [Novikov 1981]. Since $G/\ker u$ is free Abelian of rank $r={\rm rk}(u)$, $\Z[G/\ker u]_{\overline u})$ is Abelian. 
\smallskip
If $H\triangleleft_{f.i} G$ is a normal subgroup with finite index, we generalize the above definition. Consider the induced morphism
$\overline u:G/(H\cap\ker u)\to\R$ and the associated Novikov ring $\Z/[G/(H\cap\ker u)]_{\overline u}$, and define the {\it  Abelianized relative Novikov homology}
$$H_1^{ab}(G,H,u):=H_1(G;\Z[G/(H\cap \ker u)]_{\overline u}).$$

\medskip
Now we can state and prove the relation between Alexander polynomials and Abelianized relative Novikov homology. 

\medskip
\begin{proposition} 
\label{novikovandmonic} We have the equivalence

$$\Delta^H(G,u)\ {\rm is}\ u-{\rm monic} \Leftrightarrow H_1^{ab}(G,H,u)=0.$$

\end{proposition}

\medskip
\noindent{\it Proof (inspired by the referee).} Set $R=\Z[G/\ker u]$, $S=\Z[G/\ker u]_{\overline u}$, which are UFDs with $R\subset S$. 
Set

$$M_R=\Z[G/\ker u]^{G/H}\subset M_S=\Z[G/\ker u]_{\overline u}^{G/H}.$$
Then $M_R$ and $M_S$ are $G$-modules for the natural actions on $G/\ker u$ and on $G/H$, and $H_1(G;M_R)$ and $H_1(G;M_S)$ can be viewed as modules over $R$ and $S$ respectively, with $G/\ker u$ acting without permuting the factors.
Moreover, we have
$$\Delta^H(G,u)={\rm ord}_R(H_1(G;M_R))\ , \ H_1^{ab}(G,H,u)=H_1(G;M_S).$$
Recall that for an element of $S$, in particular of $R$, to be $u$-monic means to be a unit in $S$. Denote by $S^*$ the units of $S$. Using Proposition \ref{orderandannihilator}, it suffices to prove that 
$${\rm ann}_R(H_1(G;M_R))\cap S^*\ne\emptyset \Leftrightarrow H_1(G;M_S)=0.\leqno(*)$$

\smallskip
Using a presentation $\langle x_1,\cdots,x_p\mid r_1,\cdots,r_q \rangle$ be a presentation of $G$, we have
an exact complex
$$C_2=\Z[G]^q\stackrel{\times D_2}{\longrightarrow}C_1=\Z[G]^p\stackrel{\times D_1}{\longrightarrow} C_0=\Z[G] \stackrel{\varepsilon}{\longrightarrow} \Z\to 0,$$
where $D_1=\begin{pmatrix}x_1-1\\
\cdots\\
x_p-1
\end{pmatrix}$ and $\varepsilon$ is the augmentation. Thus
 $H_1(G;M_R)$ and $H_1(G;M_S)$ are the $H_1$ of the induced sequences with $\Z[G]$ replaced by 
$M_R$ and $M_S$. Thus we have (up to isomorphisms)
$$H_1(G;M_R)=\frac{\ker (\times B)}{{\rm im}(\times A)}\ , \ 
H_1(G;M_S)=\frac{\ker (\times B_S)}{{\rm im}(\times A_S)},$$
where $A\in{\rm M}_{a,b}(R)$, $B\in {\rm M}_{b,c}(R)$, and $A_S=A$, $B_S=B$ viewed as matrices with coefficients in $S$. Also, $a=q[G:H]$, $b=p[G:H]$, $c=[G:H]$ and $B=\begin{pmatrix}(\overline x_1-1){\rm I_c}\\\cdots\\ (\overline x_p-1){\rm I}_c\end{pmatrix}$.

\smallskip
The key point is that, since $u\ne 0$, there exists $i$ such that $u(x_i)\ne 0$, thus the image $\overline x_i-1\in R$ belongs to $S^*$. We can assume that $i=p$, thus 
$$B=\begin{pmatrix}B_1\\ (\overline x_p-1){\rm I}_c\end{pmatrix}\ {\rm with}\ B_1\in {\rm M}_{b-c,c}(R).$$
  Write $A=(A_1\mid A_2)$, where $A_1\in {\rm M}_{a,b-c}(R)$ and $A_2\in{\rm M}_{a,c}(R)$.
  Similarly, we have $A_S=((A_1)_S\mid (A_2)_S)$. Since $AB=0$, we have
$$A_1B_1+(\overline x_p-1)A_2=0.\leqno(*)$$

\smallskip
Proposition \ref{novikovandmonic} will result from the following lemma, a variant of \ref{deletingonecolumn}.

\medskip
\begin{lemma}
\label{deletingonecolumnbis}
{\it Up to isomorphisms, we have
$$\begin{aligned}&(\overline x_p-1){\rm coker}(\times A_1)\subset H_1(G;M_R)\subset {\rm coker}(\times A_1)\\
&H_1(G;M_S)={\rm coker}(\times (A_1)_S).
\end{aligned}$$
}

\end{lemma}

\noindent{\it Proof of Lemma \ref{deletingonecolumnbis}.}  The secund line follows from the first by replacing $R$ by $S$ and using the fact that $\overline x_p-1$ is invertible in $S$. It can also be proved directly as  Corollary \ref{deletingonecolumn}. Thus it suffices to prove the first line. 

\medskip
1) Consider the map 
$$i:w=(\overline x_p-1)v\in (\overline x_p-1)R^{b-c}\mapsto((\overline x_p-1)v,vB_1)\in \ker (\times B).$$
If $z\in R^a$, we have
$$i(w)=zA\Leftrightarrow (\overline x_p-1)v=zA_1\ {\rm and} \  vB_1=zA_2.$$
In view of $(*)$, the right hand side is equivalent to ($(\overline x_p-1)v=zA_1$). Thus $i^{-1}({\rm im}(\times A))\subset{\rm im}(\times A_1)$, thus $i$ induces an injection
$$(\overline x_p-1){\rm coker}\ (\times A_1)\to \frac{\ker (\times B)}{{\rm im}\ (\times A)}=H_1(G;M_R).$$

\medskip
2) Denoting an element of $R^b$ by $(v,w)$ with $v\in R^{b-c}$ and $w\in R^c$, consider the map
$$\pi:(v,w)\in \ker (\times B)\mapsto v\in R^{a-c}.$$
If $\pi(v,w)\in {\rm im}(\times A_1)$, i.e. there exists $z\in A^a$ such that $zA_1=v$, we have 
$$zA-(v,w)=(0,zA_2-w).$$
Since $AB=0$ and $(v,w)B=0$, this implies
$$(0,zA_2-w)B=(\overline x_p-1)(zA_2-w)=0,$$
 thus $w=zA_2$. Thus 
$$(v,w)=(zA_1,zA_2)=zA.$$
Thus $\pi^{-1}({\rm im}(\times A_1))\subset{\rm im}(\times A)$, thus $\pi$ induces an injection 
$$H_1(G;M_R)=\frac{\ker (\times B)}{{\rm im}\ (\times A)}\to {\rm coker}\ (\times A_1).$$

\medskip
\noindent{\it End of the proof of Proposition \ref{novikovandmonic}.} 
Since $\overline x_p-1\in S^*$, the lemma implies
$$\begin{aligned}
{\rm ann}_R(H_1(G;M_R))\cap S^*\ne\emptyset\ &\Leftrightarrow (\exists \lambda\in R\cap S^*) \ \lambda R^b\subset M_R^{b-c}A_1\\
&\Leftrightarrow (\exists \lambda\in R\cap S^*\ , \ X\in {\rm M}_{b-c,a}(R))\ XA_1=\lambda{\rm I}_{b-c} & (1)\\
H_1(G;M_S)=0 \ (1)&\Leftrightarrow (\exists\widehat X\in {\rm M}_{b-c,a}(S))\ \widehat X(A_1)_S={\rm I}_{b-c}. & (2)
\end{aligned}$$
Clearly, (1) $\Rightarrow$ (2). Conversely, if $\widehat X(A_1)_S={\rm I}_{b-c}$, truncating $\widehat X$ below a sufficiently high level of $u$ gives an identity with coefficents in $R$:
$$YA_1={\rm I}_{b-c}+C\ {\rm with}\ u>0 \ {\rm on}\ {\rm supp}(C).$$
Thus $\det(YA_1)\in R\cap S^*$, which implies that 
$X=YA_1\widetilde{(YA_1})^T$ (transpose of the cofactor matrix) satisfies (1).
This finishes the proof of proposition \ref{novikovandmonic}.

\section{Computations in dimension three and reduction of the main result}
\label{computationsindimensionthree}

In this section we consider the case where $G=\pi_1(M)$ where $M$ is a closed and connected three-manifold with a contractible universal covering $\widetilde M$. 

\subsection{A convenient complex for computing the homology of $G$}

Using a handle decomposition of genus $p$ [or a self-indexing Morse function with one minimum and one maximum], one can obtain $H_*(\widetilde M;\Z)$ by a complex of left modules over $\Z[G]$, of the form
$$C_*=(\Z[G]\stackrel{\times D_3}{\longrightarrow}\Z[G]^p\stackrel{\times D_2}{\longrightarrow}Z[G]^p\stackrel{\times D_1}{\longrightarrow}\Z[G])$$
with
$$D_1 =\begin{pmatrix} x_1-1\\
\cdots\cr x_p-1\\
\end{pmatrix}\ ,  \ D_3 =(y_1 -1\mid\cdots\mid y_p -1).$$
$(x_1,\cdots,x_p)$ and $(y_1,\cdots,y_p)$ are generating systems for $G$. Since $u\ne 0$, we can reorder them
so that $u(x_p)$ and $u(y_p)$ are nonzero, thus $x_p - 1$ and $y_p - 1$ are invertible in $\Z[G]_u$. Then we denote by $c$ the column $(x_i-1)_{i<m}$ and $\ell$ the row $(y_j-1)_{j<m}$, so that
$$D_1=\begin{pmatrix} c\\
 x_p-1\\
 \end{pmatrix}\ , \ D_3=(\ell \mid y_p -1).$$
Without using the contractibility of $\widetilde M$, this complex gives a free resolution of $\Z$ over $\Z[G]$ up to degree $2$, thus $$H_1(G,u)\approx\ker(\partial_1)_{\Z[G]_u} /{\rm im} (\partial_2)_{\Z[G]_u},$$ where we have changed the coefficients from $\Z[G]$ to $\Z[G]_u$. 
\smallskip
By the contractibility of $\widetilde M$, it is a complete free resolution, and $G$ is a group with $3$-dimensional Poincar\'e duality, which can be expressed as follows. Denote by $w$ the orientation morphism $G\to\{1,-1\}$, and define modified adjoint isomorphisms 
$$\lambda=\sum_g a_gg\mapsto\lambda^*=\sum_g a_g\varepsilon(g)g^{-1}\ , \ A=(a_{i,j})^*=(a_{j,i}^*).$$
Then Poincar\'e duality can be expressed by the fact that $C_*$ is quasi-isomorphic to the complex $(C^*_i = C_{3-i},\times D_{4-i}^*)$.

\medskip
Let us write $D_2 = \begin{pmatrix} A&C\\
 L&a\\ \end{pmatrix}$  where A $ \in {\rm M}_{p-1}(\Z[G])$, $L  \in M_{1,p-1}(\Z[G])$, $C  \in {\rm M}_{p-1,1}(\Z[G])$ and $a  \in \Z[G]$. Note that $D_2^{(p)}=(A\mid C)$.
\smallskip
Since $\partial_1\circ\partial_2 =0$ and $\partial_2\circ\partial_3 =0$, we have $D_2D_1 =0$ and $D_3D_2 =0$. Working over $ \Z[G]_u$, we
obtain
$$\begin{array}{ll}
C&=Ac(1-x_p)^{-1}\ , \ a=Lc(1-x_p)^{-1}\\
L&=(1-y_p)^{-1}\ell A\\
a&=(1-y_p)^{-1}\ell C=(1-y_p)^{-1}\ell Ac(1-x_p)^{-1}.
\end{array}$$
Thus $(D_3^*,D_2^*,D_1^*)$ has ``the same shape'' as $(D_1,D_2,D_3)$ in the following sense: it  is obtained from $(D_1,D_2,D_3)$  by replacing $(x_i,y_i,A,C,L,a)$ by $(y_i^{-1},x_i^{-1},A^*,L^*,C^*,a^*)$, and one has $u(y_p)\ne0$.

\smallskip
These computations have the following consequence.

\begin{proposition} 
\label{equivalencefiberinginvertibility}  Let $u\in H^1(M;\R)\setminus\{0\}$. The two following properties are equivalent:

\begin{enumerate}
\label{bierietal}

\smallskip
\item $u$ is fibered.

\smallskip
\item The matrix $A\in{\rm M}_{p-1}(\Z[G])$ becomes invertible in ${\rm M}_{p-1}(\Z[G]_{u})$.
\end{enumerate}
\end{proposition}

\medskip
\noindent{\it Proof.} By \cite{Bieri-Neumann-Strebel 1987}, (1) is equivalent to ($H_1(G,u)=0$) thus to the exactness of $C_*\otimes_{\Z G}\Z[G]_u$ in degree $1$. Since $u(x_p)\ne0$, $\overline x_p-1$ is a unit of $\Z[G]_u$. By \ref{deletingonecolumn}, this is equivalent to the left invertibility of $\begin{pmatrix} A\\ L \end{pmatrix}$ with $\Z[G]_u$-coefficients. Since $L $ is of the form $\lambda A$, this left invertibilty is equivalent to that of $A$ in ${\rm M}_{p-1}(\Z[G]_u)$.

\smallskip
Since $C^*$ is quasi-isomorphic to $C_*$, (1) is equivalent to the exactness of $C^*\otimes_{\Z G}\Z[G]_u$ in degree $1$. Since  $(D_3^*,D_2^*,D_1^*)$ has the same shape as $(D_1,D_2,D_3)$ in the sense explained above, we can apply the same argument to prove that (1) is equivalent to the left invertibility of $A^*$ in ${\rm M}_{p-1}(\Z[G]_{-u})$,  i.e.  the right-invertibility of $A$ in ${\rm M}_{p-1}(\Z[G]_{u})$. This proves Proposition \ref{bierietal}.

\medskip

\noindent{\it Remark.}   In \cite{Sikorav 2017} it was established that $\Z[G]_u$ is always {\it stably finite}: a matrix $A  \in {\rm M}_n(\Z[G]_u)$ is invertible if and only if it is left invertible. This is a well-known result of Kaplansky for $\Z[G]$, which was proved by \cite{Kochloukova 2006} for $\Z[G]_u$ when rk$(u) = 1$.  With Poincar\'e duality, this allows to prove ($H_1(G,u)=0\Rightarrow H_1(G,-u)=0$) without using the results of Stallings and Thurston.

\medskip
Proposition \ref{equivalencefiberinginvertibility} reduces the proof of the main result  to the following result.

\begin{theorem}  
\label{mixedfinitedetectability}
Assume that $G$ is finitely generated and VRTFN. Let $p$ be a prime and $n\in\N$. Let $A\in{\rm M}_n(\Z[G])$
be such that for every $H\triangleleft_{f.i.} G$ its image in ${\rm M}_n(\Z[G/H \cap  \ker u]_{\overline u})$ is invertible.
\smallskip
Then $A$ is invertible in ${\rm M}_n(\Z[G]_u)$.
\end{theorem}

\begin{remarks}
\label{finiteindex}
{\rm \medskip
\begin{enumerate} 

\item Note that we state the result for $A\in{\rm M}_n(\Z[G])$, not in ${\rm M}_n(\Z[G]_u)$. Presumably, the result would remain true, but it is not needed and I have not been able to prove it.

\smallskip
\item
The validity of Theorem \ref{mixedfinitedetectability} for a finite index subgroup $G_0\subset G$ implies its validity for $G$: this follows from the fact that an $n$-matrix $A$ over $\Z[G]_u$ can be represented by a $(n.[G\!:\!G_0])$-matrix $\widetilde A$ over $\Z[G_0]_{u_{|G_0}}$, and that the invertibility of $A$ is equivalent to the bijectivity of the left and right multiplications by $A$, thus to the invertibility of $\widetilde A$ (and similarly for the finite invertibility).

\smallskip
Thus it suffices to prove Theorem \ref{mixedfinitedetectability} when $G$ is finitely generated and RTFN.
\end{enumerate}

}
\end{remarks}

\bigskip
\section{Finitely detectable units and full left ideals in group rings}
\label{finitelydetectableunitsandfullleftideals}

\subsection{Definitions.}  \begin{enumerate}

\smallskip
\item A matrix $A\in {\rm M}_n(\Z[G])$ is {\it finitely invertible} if its image in every quotient ${\rm M}_n(\Z[G/H])$  for $H\triangleleft_{f.i.} G$ is invertible.  The ring $ {\rm M}_n(\Z[G])$ has {\it finitely detectable units} if every finitely invertible matrix is invertible.

\smallskip
\item A left ideal $I\subset \Z[G]$ is {\it finitely full} if the natural projection $I_H\subset \Z[G/H]$ is equal to $\Z[G/H]$ for every $H\triangleleft_{f.i.} G$.  The ring $\Z[G]$ has {\it finitely detectable full left ideals} if every left ideal which is finitely full is equal to $\Z[G]$.

\end{enumerate}

\begin{remark} {\rm Since $\Z[G]$ is anti-isomorphic to itself via $\sum a_gg\mapsto\sum a_gg^{-1}$, if $\Z[G]$ has finitely detectable full left ideals, it also has detectable full right ideals.
}
\end{remark}

\begin{proposition} Assume that $\Z[G]$ has finitely detectable full left ideals.

\begin{enumerate}

\smallskip
\item Every left $\Z[G]$-submodule ${\mathcal M}\subset \Z[G]^n$ which projects onto $\Z[G/H]^n$ for every $H\triangleleft_{f.i.} G$ is equal to $\Z[G]^n$.

\smallskip
\item For every $n\in\N^*$, ${\rm M}_n(\Z[G])$ has finitely detectable units.
\end{enumerate}

\end{proposition}

\medskip
\noindent{\it Proof. } (1) For $n = 1$, it is the hypothesis. Assume that $n > 1$ and the result is true for $n-1$. Consider the set $I$ of $\lambda  \in \Z[G]$ such that there exists $\lambda_1,\cdots,\lambda_{n-1}  \in \Z[G]$ with $(\lambda_1,\cdots,\lambda_{n-1},\lambda)  \in {\mathcal M}$. It is a left ideal, which projects onto every quotient $\Z[G/H]$ with $H\triangleleft_{f.i.} G$. Thus $I = \Z[G]$,  i.e. ${\mathcal M}$ contains an element $x=(\lambda_1,\cdots,\lambda_{n-1},1)$.

\smallskip
One has a direct sum decomposition 
$$\Z[G]^n = \Z[G]^{n-1} \oplus \Z[G] x.$$ Subtracting $\mu_nx$ from every element $(\mu_1,\cdots,\mu_n)  \in {\mathcal M}$, one sees that 
$${\mathcal M} = ({\mathcal M} \cap \Z[G]^{n-1})\oplus \Z[G]x.$$
It suffices to prove
that ${\mathcal M} \cap  \Z[G]^{n-1} = \Z[G]^{n-1}$. Clearly, ${\mathcal M} \cap  \Z[G]^{n-1}$ is a left submodule which projects onto every
$(\Z[G/H])^{n-1}$, $H\triangleleft_{f.i.} G$. Thus the proposition follows from the induction hypothesis.

\medskip
(2) If $A  \in {\rm M}_n(\Z[G])$ is finitely invertible, the left submodule $\Z[G]^nA\subset \Z[G]^n$ is finitely full. Thus $\Z[G]^nA= \Z[G]^n$, i.e.
$A$ is left invertible. Similarly, $A\Z[G]^n= \Z[G]^n$, thus $A$ is right invertible. Thus $A$ is a unit.$\qed$

\bigskip
\section{Facts on nilpotent groups and their group rings}
\label{factsonnilpotent}

We collect here a few facts that we will use about (mostly finitely generated) nilpotent groups and their group rings. 

\medskip

\begin{enumerate}

\medskip

\item
If $G$ is nilpotent and finitely generated, $G$ is polycyclic.
 (\cite{Kargapolov-Merzliakov}, Theorem 17.2.2 p.119; \cite{Robinson}, 5.2.17 p.13.)
\medskip
\item If $G$ is nilpotent and finitely generated, it has a torsion-free subgroup of finite index.\break (\cite{Kargapolov-Merzliakov}, Theorem 17.2.2 p.119.)

\medskip
\item If $G$ is nilpotent and finitely generated, it is residually finite. (\cite{Kargapolov-Merzliakov}, Exercise 17.2.8 p.124 (follows from (2) and the Mal'cev embedding of a finitely generated torsion-free nilpotent group in some ${\rm SL}(n,\Z)$, Theorem 17.2.8 p.120.)

\medskip
\item If $G$ is nilpotent and finitely generated, every subgroup of $G$ is finitely generated (\cite{Robinson}, 5.2.18 p.137.)

\medskip
\item  If $G$ is any group and $(\gamma_n(G))$ its lower central series ($\gamma_1(G)=G$, $\gamma_{n+1}(G)=[G,\gamma_n(G)]$), the set 
$$G_n:=\sqrt{\gamma_n(G)}:=\{g\in G\mid (\exists k\in\N^*)\ g^k\in\gamma_n(G)\}$$ is a normal subgroup, moreover $G/G_n$ is torsion-free and $[G_n,G_m]\subset G_{n+m}$. (\cite{Passman}, Lemma 1.8 p.473.)
Note that the sequence $(G_n)$ is finite iff $G$ is nilpotent and torsion-free.

\medskip
\item If $G$ is nilpotent and torsion-free, it is {\rm orderable}, i.e. it has a total order such that $x\le y\Rightarrow xz\le yz$ and $zx\le zy$. {\rm (See also Section \ref{section:malcevneumanncompletion}.)} (\cite{Passman}, Lemma 1.6 p.587.) 

\medskip \item If $G$ is nilpotent and torsion-free, $\Z[G]$ is a domain. (Easy consequence of (6).)

\medskip
\item If $G$ is polycyclic (in particular, nilpotent and finitely generated), $\Z[G]$ is left (and right) Noetherian. (\cite{Hall 1954}; \cite{Passman}, Corollary 2.8 p.425.)

\medskip
\item If $G$ is nilpotent and finitely generated, a left $\Z[G]$-module which is simple (or irreducible) is finite. Equivalently, if $I$ is a maximal left ideal of $\Z[G]$, $\Z[G]/I$ is finite.  (\cite{Hall 1959}, Lemma 2 and Theorem  3.1;  \cite{Passman}, Corollaries 2.9 and 2.10 p.544.) 

\medskip
\item If $G$ is nilpotent and torsion-free, $\Z[G]$ is a {\it UFD in the noncommutative sense} \cite{Chatters 1984}. This means that every nonzero element is a product of irreducible elements, the decomposition being unique up to order and multiplication by units. Moreover, an irreducible element $\lambda$ satisfies $\lambda\Z[G]=\Z[G]\lambda$.

\medskip

\item If $G$ is nilpotent and torsion-free, $\Z[G]$ is contained in a division ring $D$ which is a {\rm classical ring of quotients on the right and on the left}:
$$D= \{xy^{-1}\mid x\in  \Z[G]\ , \ y\in \Z[G]\setminus\{0\}\}= \{y^{-1}x\mid x\in  \Z[G]\ , \ y\in  \Z[G]\setminus\{0\}\}.$$
(\cite{Goldie 1958}, Theorem 1 and \cite{Lam}, Corollary 10.23 p.304: A right Noetherian domain has a classical ring of quotients).
\end{enumerate}

The last statement has the following consequence.

\begin{corollary} (\cite{Lam} p.301)
\label{commondenominator} If $G$ is nilpotent and torsion-free and $E\subset \Z[G]$ is finite, its elements can be reduced to a common denominator: there exists $x\in \Z[G]\setminus\{0\}$ such that $E\subset \Z[G]x^{-1}$.
\end{corollary}

\bigskip
\section{Theorem 7.3 for nilpotent groups}
\label{section:theorem7.3fornilpotent}

\medskip
In this section, $G$ is a finitely generated nilpotent group. 
We first prove the finite detectability of full ideals, then the result in the title.

\begin{proposition} 
\label{finitedetectabilityfornilpotent} Let $G$ be a finitely generated nilpotent group. Then $\Z[G]$ has finitely detectable full left ideals. 

\end{proposition}

\medskip
\noindent{\it Proof.} (1) We argue by contradiction, thus we assume that $I$ is a finitely full left ideal in $\Z[G]$ which is not full. Then $I$ is contained in a maximal left ideal $I_1$, without the axiom of choice since $\Z[G]$ is Noetherian.
Then $I_1$ is again a finitely full left ideal in $\Z[G]$ which is not full, thus we can assume that $I$ is maximal.

\smallskip
 Thus  ${\mathcal M} := \Z[G]/I$ is a simple $\Z[G]$-module, and by \cite{Hall 1959}, ${\mathcal M}$ is finite.Thus
$$H:=\ker(G\to{\rm Aut}({\mathcal M}))$$
has finite index. Thus ${\mathcal M}$ is isomorphic to a quotient of $\Z[G/H]/I_H$ with $H\triangleleft_{f.i.}G$ and $I_H$ the image of $I$. By hypothesis, $I_H=\Z[G/H]$, thus ${\mathcal M}=0$, contradiction.

\medskip
Now we prove Theorem 7.3 for nilpotent finitely generated groups.

\begin{proposition}
\label{theorem7.3fornilpotent} Let $G$ be a nilpotent and finitely generated group, and let $A\in{\rm M}_m(\Z[G])$
be such that, for every $H\triangleleft_{f.i.} G$, the image $A_{H,u}$ of $A$ in ${\rm M}_m(\Z[G/H \cap  \ker u]_{\overline u})$ is invertible. Then $A$ is invertible in ${\rm M}_m(\Z[G]_u)$.

\end{proposition}

\medskip
\noindent{\it Proof.} By Remark \ref{finiteindex} (3), we can assume that $G$ is torsion-free, thus
$\Z[G]$ is contained in a division ring $D$ which is a classical ring of fractions on the right and on the left.

\smallskip
If $u$ is injective, the result is obvious. In general, we make an induction over the polycyclic length of $G$, i.e. in the torsion-free case the length $n$ of any subnormal 
sequence of subgroups
$$G=G_0>G_1>\cdots>G_n>G_{n+1}=\{1\}\ , \ G_i/G_{i+1}\approx\Z.$$
By the Schreier refinement theorem (\cite{Kargapolov-Merzliakov}, 4.4.4), this length is independent of the sequence.
\smallskip
We can assume that $\ker u\ne\{1\}$ and that the Proposition is already known when the Hirsch length is smaller than 
that of $G$. Then $u$ is not injective on the center $C(G)$, otherwise we would have $C(G)\cap [G,G]=\{1\}$
thus $[G,G]=\{1\}$ and $G$ would be Abelian, giving a contradiction. 

\smallskip
Moreover, every element of $G$ which has a nontrivial power in $C(G)$ is already in $G$ (\cite{Kargapolov-Merzliakov}, Exercise 16.2.9).  Thus we can find $z\in C(G)\cap\ker u$ such that $\Gamma:=G/\langle z\rangle$ has no torsion. Then the Hirsch length of $\Gamma$ is smaller than that of $G$.

\begin{lemma}
\label{AB=xI}
There is an identity $AB=x{\rm I}_m$ in ${\rm M}_m(\Z[G]_u)$, with $x\in\Z[G]_u\setminus\{0\}$.

\end{lemma}

\noindent{\it Proof.}   The right multiplication by $A$ is injective on $(\Z[G]_u)^n$: if $LA=0$, we obtain $L_{H,u}A_{H,u}=0$ for every $H\triangleleft_{f.i.} G$, where $L_{H,u}$ is the image of $L$ in $(\Z[G/(H\cap\ker u)]_{\overline u})^m$. 
Since by hypothesis $A_{H,u}$ is invertible, $L_{H,u}=0$. Since this is true for all $H$ and $G$ is residually finite, $L=0$. 

\smallskip
Thus $A$ has an inverse in ${\rm M}_m(D)$, which by Corollary \ref{commondenominator} is of the form $A^{-1}=Bx^{-1}$ with $B\in{\rm M}_m(\Z[G])$, $x\in \Z[G]\setminus\{0\}$. This proves Lemma \ref{AB=xI}.

\medskip
To prove Proposition \ref{theorem7.3fornilpotent}, it suffices to prove that $x$ divides $B$  on the right, i.e. $B=B_1x$ in ${\rm M}_n(\Z[G]_u)$, thus $(AB_1)x=x{\rm I}_n$, and since $x\ne 0$ and $\Z[G]_u$ is a domain, this implies $AB_1={\rm I}_n$. A similar argument with ``left'' and ``right'' exchanged proves that $A$ is left-invertible, thus invertible.

\begin{lemma} 
\label{xdividesB}
For $n\in\N^*$, let $x_n,B_n$ be the images of $x,B$ in $\Z[G/\langle z^n\rangle]_{\overline u})$ and ${\rm M}(m,\Z[G/\langle z^n\rangle]_{\overline u}))$. Then
$x_n$ divides $B_n$ on the right.
\end{lemma}

\medskip
\noindent{\it Proof.} The image $A_n$ of $A$ in ${\rm M}_m(\Z[G/\langle z^n\rangle])$ gives rise to a matrix $\widetilde A_n\in{\rm M}_{mn}(\Z[\Gamma])$ whose images in every ${\rm M}_{mn}(\Z[\Gamma/H\cap\ker\overline u]_{\overline u})$ is invertible, and $\widetilde A_n$ is invertible in ${\rm M}_{mn}(\Z[\Gamma]_{\overline u})$ if and only if $A_n$ is invertible in ${\rm M}_m(\Z[G/\langle z^n\rangle]_{\overline u})$.

\smallskip
By the induction hypothesis, $\widetilde A_n$ is invertible in ${\rm M}_{mn}(\Z[\Gamma]_{\overline u})$ thus $A_n$ is invertible in ${\rm M}_m(\Z[G/\langle z^n\rangle]_{\overline u})$. Denote its inverse by $A_n^{-1}$ and multiply the identity $A_nB_n=x_n{\rm I}_m$  
on the left by $A_n^{-1}$, we obtain $B_n=A_n^{-1}x_n$, which proves Lemma \ref{xdividesB}.

\medskip
We shall need the two following objects.

\begin{enumerate}

\medskip
\item
For $\lambda=\displaystyle\sum_{g\in G} a_gg\in\Z[G]_u\setminus\{0\}$, define $\mu=\min(u_{|{\rm supp}(z)})$ and 
$$\widetilde m_u(\lambda)=\displaystyle\sum_{g\mid u(g)=\min(u_{|{\rm supp}(\lambda)})}a_gg\in\Z[G].$$
We have $\widetilde m_u(\lambda)=gm_u(\lambda)$ with $g\in G$ and $m_u(\lambda)\in\Z[G/\ker u]$, where $m_u(\lambda)$ is defined up to multiplication by an element of $\ker u$. We call $m_u(\lambda)$ the {\it $u$-minimal part of $\lambda$}.

\medskip
\item
Let $\zeta_n\in\C$ be a primitive $n$-root of unity, and let $\Phi_n\in\Z[t]$ be its minimal polynomial (the $n$-th cyclotomic polynomial). The rings $\Z[G]$ and $\Z[G]_{\overline u}$ can be factored by the ideal generated by $\Phi_n(z)$, to give quotients of $\Z[G/\langle z^n\rangle]$, and $\Z[G/\langle z^n\rangle]_{\overline u}$. The quotients may be expressed as twisted rings
$$\begin{aligned}\Z[G]/(\Phi_n(z))&=\Z[\zeta_n][\Gamma]\\
\Z[G]_{\overline u}/(\Phi_n(z))&=\Z[\zeta_n][\Gamma]_{\overline u}.
\end{aligned}$$
Since $\Z[\zeta_n]$ is a domain, these rings are also domains. 

\end{enumerate}

\medskip
Let $y$ be a coefficient of $B$. Denote by $\overline x_n$ and $\overline y_n$ the images of $x$ and $y$ in $\Z[\zeta_n][\Gamma]_{\overline u}$. By Lemma \ref{xdividesB}, $\overline x_n$ divides $\overline y_n$ since they are also images of $x_n,y_n\in\Z[G/\langle z^n\rangle]_{\overline u}$. Moreover, $\overline x_n$ and $\overline y_n$ have $\overline u$-minimal parts 
$$m_{\overline u}(\overline x_n),m_{\overline u}(\overline y_n)\in\Z[\zeta_n][\ker\overline u],$$ defined up to multiplication by an element of $\pm \ker\overline u$. 

\smallskip
Since $m_u(x)\ne0$, for $n>\!\!>1$ its image $\overline {m_u(x)}_n$ is nonzero, thus equal to
$m_{\overline u}(\overline x_n)$. And since $\overline x_n$ divides $\overline y_n$ and $\Z[\ker\overline u]$ is a domain, this implies that $\overline {m_u(x)}_n$ 
divides $\overline{m_u(y)_n}$ for $n>\!\!>1$, equivalently that $m_u(x)$ divides $m_u(y)$ modulo $\Phi_n(z)$ for $n>\!\!>1$. 

\medskip
 To finish the proof of Proposition \ref{theorem7.3fornilpotent}, it suffices to prove the following lemma.

\begin{lemma} \label{divides}
\

\begin{enumerate} 

\medskip
\item In $P,Q\in\Z[t]$ and $P(\zeta_n)$ divides $Q(\zeta_n)$ in $\Z[\zeta_n]$ for $n>\!\!>1$, then $P$ divides $Q$ in $\Z[t]$.

\medskip

\item If $\lambda,\mu\in\Z[\ker u]$ and $\lambda$ divides $\mu$ modulo $\Phi_n(z)$ for $n>\!\!>1$, then $\lambda$ divides $\mu$ in $\Z[\ker u]$.

\medskip
\item If $x,y\in\Z[G]_u$ and $m_u(x)$ divides $m_u(y)$ modulo $\Phi_n(z)$ for $n>\!\!>1$, then $x$ divides $y$ in $\Z[G]_u$.

\end{enumerate}
\end{lemma}

\smallskip
\noindent Indeed, modulo the lemma we have $B=B_1x$ with $B_1\in{\rm M}_n(\Z[G]_u)$, thus $(AB_1)x=x{\rm I}_n$, and since $x\ne0$ and $\Z[G]_u$ is a domain, this implies $AB_1={\rm I}_n$, thus $A$ is right invertible. A similar argument as above with ``right" and ``left" exchanged proves that $A$ is left-invertible, thus $A$ is invertible.

\medskip
\noindent{\it Proof.}  

\medskip
(1) Since $\Z[t]$ is a UFD, one can reduce to the case when $P$ is irreducible. The resultants 
${\rm res}(P,\Phi_n)$ and ${\rm res}(P,Q)$ satisfy
$${\rm res}(P,\Phi_n)= \pm\prod_{\zeta\in \Phi_n^{-1}(0)}P(\zeta_n)\ , \ {\rm res}(P,Q)=\pm\prod_{\zeta\in \Phi_n^{-1}(0)}P(\zeta_n).$$
Since $P(\zeta_n)$ divides $Q(\zeta_n)$ in $\Z[\zeta_n]$ for $n>\!\!>1$, this implies$$(\forall n>\!\!>1)\  {\rm res}(P,\Phi_n)\  {\rm divides}\ {\rm res}(P,Q)\ {\rm in}\ \Z.$$
We want to prove that ${\rm res}(P,Q)=0$. It suffices to prove that ${\rm res}(P,\Phi_n)$ takes infinitely many values.

\smallskip
If $p$ is a prime number, we have 
$${\rm res}(P,\Phi_p)={\rm res}(P,t^{p-1}+\cdots+1)=\prod_{\alpha\in P^{-1}(0)}(\alpha^{p-1}+\cdots+1).$$
We distinguish three cases:

\medskip
$\bullet$ The  zeros $\alpha_1,\cdots,\alpha_d$ of $P$ are not algebraic units. Then for some non-Archimedean absolute value $|.|_v$ on $\Q(\alpha_1,\cdots,\alpha_d)$ we have  $|\alpha_1|_v=\cdots=|\alpha_d|_v\ne 1$. Replacing $P$ by $t^dP(t^{-1})$, we can assume that $|\alpha_1|_v>1$, thus as $p\to\infty$ the formula for ${\rm res}(P,\Phi_p)$ implies that when $p\to\infty$ we have
$$|{\rm res}(P,\Phi_p)|_v\sim C|\alpha_1|^{dp}\ , \ C>0.$$
Thus ${\rm res}(P,\Phi_p)$ takes  infinitely many values.

\medskip
$\bullet$  The zeros of $P$ are algebraic units but not roots of unity. Then at least one has modulus $1$. Say that $|\alpha_1|,\cdots,|\alpha_k|>1\ge|\alpha_{k+1},\cdots,\alpha_d$. Then the formula for ${\rm res}(P,\Phi_p)$ implies
$$|{\rm res}(P,\Phi_p)|\sim C|\alpha_1|^{kp}\ , \ C>0.$$
Again,  ${\rm res}(P,\Phi_p)$ takes infinitely many values.

\medskip
$\bullet$ The zeros of $P$ are roots of unity, i.e. $P=\pm\Phi_k$ for some $k$. Then if $p$ is a large prime, it does not divides $k$, thus
$$\Phi_{kp}=\sum_{i=0}^{p-1}t^{ki},$$
which implies
$${\rm res}(\Phi_k,\Phi_{kp})=\prod_{\zeta\in \Phi_k^{-1}(\{0\})}\Phi_{kp}(\zeta)=p^{\varphi(k)}.$$
Thus ${\rm res}(P,\Phi_{kp})$ takes infinitely many values.

\medskip
 (2) Since $\Z[\ker u]$ is a UFD in the non commutative sense, one can reduce to the case when $\lambda$ is irreducible. Then $R=\Z[G]/\lambda\Z[G]$ is a ring and a domain. It is also Noetherian on the right and on the left. Let $\overline z,\overline\mu$ be the images of $z,\mu$ in $R$, by hypothesis we have that $\Phi_n(\overline z)$ divides $\overline\mu$ for $n$ large. If $\overline\mu=0$ we are done. Assume now that $\overline\mu\ne0$.

\smallskip
If $p,p'$ are distinct prime numbers, $\Phi_p$ and $\Phi_{p'}$ satisfy an identity $U\Phi_p+V\Phi_{p'}=1$ in $\Z[t]$ (more generally this is true for $1+t+\cdots+t^{r-1}$ and $1+t+\cdots+t^{s-1}$ if $r\wedge s=1$, using the Euclidean algorithm). Thus if $(p_n)$ is the sequence of prime numbers, there exist $n_0$ such that $\overline\mu$ is divisible by $\Phi_{p_{n_0}}(z)\cdots\Phi_{p_n}(z)$ for all $n\ge n_0$. Let $q_n$ be the quotient, and $I_n=q_nR$. By Noetherianity, $I_n$ is constant  for $n>\!\!>1$.

\smallskip
Thus the image of $\Phi_p(z)$ in $\Z[G]/\lambda\Z[G]$  is a unit for $p>\!\!\>1$. Equivalently, the image of $\lambda$ is a unit in $\Z[G]/(\Phi_p(z))\approx\Z[\zeta_p][\Gamma]$. Up to multiplication by a unit, we can write
$$\lambda=P(z)+\lambda_+,$$
where $P\in\Z[t]\setminus\{0\}$ and every element of the projection of ${\rm supp}(\lambda_+)$ in $G/\langle z\rangle=\Gamma$ is $>1$.
Since  the image of $P(z)+\lambda_+$ is invertible in $\Z[\zeta_n][\Gamma]$, the image of $P$ being nonzero (for $p>\!\!>1$), and since $\Z[\zeta_n]$ is a domain,  the image of $\lambda_+$ vanishes, i.e. $\lambda_+$ is divisible by $\Phi_p(z)$ for $p>\!\!>1$, thus $\lambda_+=0$.

\smallskip
Thus $\lambda=P(z)$, thus by projections on $\Z[z,z^{-1}]g$ for each $g\in G$, we can assume that $\mu=Q(z)$ with $Q\in\Z[t,t^{-1}]$. And up to multiplication by a unit, that $Q\in\Z[t]$.  Finally,  the fact that $P(z)$ divides $Q(z)$ modulo $\Phi_n(z)$ for $n>\!\!>1$ is equivalent to: $P(\zeta_n)$ divides $Q(\zeta_n)$ in $\Z[\zeta_n]$ for $n>\!\!>1$. By (1), $P$ divides $Q$, which proves (2).

\medskip
(3)  We define a ``division algorithm by increasing value of $u$" (analogous to the division algorithm in $\Q[[t]]$). Since in 2) we can replace $\mu$ by $\mu-\alpha\lambda$, we have that $m_u(x)$ divides $m_u(y-\alpha x)$ for every $\alpha\in\Z[\ker u]$. This allows to define by induction a sequence $(x_n)$ in $\Z[G]_u$ by 
$x_0=x$, $x_1=y$  and
$$(\forall n\ge2)\ \ x_n=x_{n-1}-\widetilde m_u(x_{n-1})\widetilde m_u(x)^{-1}.$$
 By construction, we have 
$$y=(\widetilde m_u(x_1)\widetilde m_u(x)^{-1}+\cdots+\widetilde m_u(x_{n-1})\widetilde m_u(x)^{-1})x+x_n\leqno(*)$$
Define $S={\rm supp}(y)$ and
$$T=\{g^{-1}h\mid h\in  {\rm supp}(\widetilde m_u(x))\ , \ g\in{\rm supp}(t_u(x))\},$$
which is a subset of $[a,+\infty[$ for some $a>0$. For $n\ge2$, we have
$${\rm supp}(x_n)\subset{\rm supp}(x_{n-1})\cup T{\rm supp}(x_{n-1}).$$
Since ${\rm supp}(x_1)=S$, this implies 
$${\rm supp}(x_n)\subset E=\displaystyle \bigcup_{i=0}^\infty ST^i\ , \ T^i=\{t_1\cdots t_i\mid t_1,\cdots,t_i\in T\}.$$
Let $\mu_n=\min(u_{|{\rm supp}(x_n)})$. By construction, $\mu_n$ is increasing and belongs to $u(E)$. Since $u\ge a>0$ on $T$, $E$ has only a finite number of elements in $\{u\le C\}$ for every $C\in\R$, thus $\mu_n\to+\infty$. Thus $\displaystyle\sum_{n=1}^\infty \widetilde m_u(x_n)\widetilde m_u(x)^{-1}$ is a well-defined element $\alpha\in \Z[G]_u$, and by $(*)$ we have $y=\alpha x$. This finishes the proof of Lemma \ref{divides} and thus of  Proposition \ref{theorem7.3fornilpotent}.

\section{Mal'cev-Neumann completion of $\Z[G]$}
\label{section:malcevneumanncompletion}

Here we assume that $G$ is residually torsion-free nilpotent (RTFN), i.e. there exists a series of normal subgroups of $G=G_0>G_1>\cdots>G_n$, such that $G/G_n$ is torsion-free nilpotent and $\displaystyle\bigcap_{n\in\N}G_n=\{1\}$.
We also require that $G$ be finitely generated (countable would suffice).

\medskip
\subsection{Order on $G$.} 
Following \cite{Eizenbud-Lichtman 1987}, we define an order on $G$ as follows.
First, one defines 
$$G_n:=\sqrt{\gamma_n(G)}=\{x\in G\mid(\exists m>0)\ x^m\in\gamma_n(G)\}.$$
As we recalled, since $G$ is nilpotent they are subgroups. Clearly, they are normal in $G$, and $G/G_n$ is torsion-free. Moreover, one has $[G_n,G_n]\subset G_{2n+1}\subset G_{n+1}$ by \cite{Passman}, Lemma 1.8 p.473. Thus $G_n/G_{n+1}$ is torsion-free Abelian. It is also finitely generated since it is contained in $G/G_n$ which is nilpotent and finitely generated.

\smallskip

One orders arbitrarily each $G_n/G_{n+1}$. Then one defines $x\in G_0$ to be positive if and only if, for the unique $n$ such that $x\in G_n\setminus G_{n+1}$, one has $xG_{n+1}>1$ in $G_n/G_{n+1}$.

In other words, an element $x\in G$ is $>1$ if and only if its first nontrivial image in a subquotient $G_{n-1}/G_n$ is $>1$.  It is clear that $G^+$ is indeed the positive cone of an order on $G$.

\subsection{Mal'cev-Neumann completion, comparison with Novikov.} 
\label{malceneumanncomparisonwithnovikov}
We recall a celebrated result of A.I. Mal'cev and B.H. Neumann: if $G$ is a bi-invariantly ordered group, the formal series
$$\Q\langle G \rangle  := \{\lambda  \in \Q[[G]]\mid {\rm supp}(\lambda)\ \hbox{\rm is well-ordered}\}$$
form a division ring (or skew field) for the natural operations, containing $\Q[G]$ as a subring. (Actually, one can replace $\Q$ by any field, or even any division ring).

\medskip
In presence of a nonzero morphism $u:G\to\R$, following \cite{Kielak 2020}, we shall require the order to be compatible with $u$ in the sense that ($u(x)>0\Rightarrow x>1$). This is possible by changing the definition of $(G_n)$, setting

\begin{itemize}

\item $G_1^{new}=\ker u$

\item $G_n^{new}=G_{n-1}$ if $n\ge2$.

\end{itemize}
and defining the order on $G/G_1$ by embedding it in $\R$ via $\overline u$ induced by $u$. The interest of this is that we then have
$$\Z[G]_u \subset \Q\langle G \rangle.$$

\subsection{Subfield with controlled coefficients.} We shall work mostly in a subfield of $\Z\langle G\rangle$, introduced by \cite{Eizenbud-Lichtman 1987}, which contains $\Z[G]$: by definition,
$$S(G):=\bigcup_{n=1}^\infty \Q[G_n]\langle G/G_n\rangle,$$
where
 $$\Q[G_n]\langle G/G_n\rangle:=\{\lambda\in \Q[[G]]\mid \lambda=\sum_{t\in T}\lambda_t t\ , \ \lambda_t\in \Q[G_n]\ , \ t\in T\},$$
where $T$ is a well-ordered subset of some transversal $T_n$ of $G_n$ in $G$ (clearly, this does not depend on the choice of $T_n$). 
This is clearly a subring of $\Q\langle G\rangle$, which contains  $\Q[G]$.

\begin{proposition} 
\label{controlledsupport} (\cite{Eizenbud-Lichtman 1987}, Proposition 4.3) $S(G)$ is a subfield of $\Q\langle G\rangle$.
\end{proposition}
We provide a proof of this proposition since that of Eizenbud-Lichtman is incorrect (I thank Andrei Jaikin (personal communication, January 2020) for alerting me about this).

\smallskip
 We can choose the transversals so that $T_n\subset T_{n+1}$. We have to prove that every nonzero $\lambda\in \Q[G_n]\langle G/G_n\rangle$ is invertible in $\Q[G_m]\langle G/G_m\rangle$ for $m$ large enough. We can assume that $\displaystyle\lambda=\sum_{t\in T}\lambda_tt$, with $T$ a well-ordered subset of $T_n$, $\min T=1$, $\lambda_t\in\Q[G_n]$ and $\lambda_1\ne0$. And also that 
 $$\lambda_1=1+a_1g_1+\cdots+a_kg_k\ , \ a_i\in\Q\ , \ g_i>1.$$
Let $m\ge n$ be large enough that none of the $g_i$ is in $G_m$. Thus $g_i=\gamma_i t_i$ with $\gamma_i\in G_m$ and $t_i\in T_m\cap G^+$. Thus
 
 $$\lambda=1+\sum_{i=1}^ka_i\gamma_i t_i+\sum_{t\in T\setminus\{t_0\}}\lambda_tt.$$
Moreover, $\lambda_tt\in \Q[G_m]F_t$ where $F_t$ is a finite subset of $T_m$, with $t<t'\Rightarrow F_t<F_{t'}$ (every element of $F_t$ is less than every element of $F_{t'}$). Since $T$ is well-ordered, $\widehat T:=\{t_1,\cdots,t_m\}\cup\displaystyle\bigcup_{t\in T} F_t$ is a well-ordered subset of $T_m\cap G^+$, and we have 
$$\lambda=1+\sum_{t\in\widehat T}\mu_tt\ , \ \mu_t\in \Q[G_m].$$
Thus in $\Q\langle G\rangle$ we have
$$\lambda^{-1}=1-\sum_{t\in\widehat T}\mu_tt+\cdots+(-1)^r\sum_{t_1,\cdots,t_r\in\widehat T}\mu_{t_1}t_1\cdots\mu_{t_r}t_r+\cdots.$$
Each product $\mu_{t_1}t_1\cdots\mu_{t_r}t_r$ can be rewritten $\nu_tt$ with $\nu_t\in \Q[G_m]$ and $t\in (\widehat T)^+$ (a positive word in $\widehat T$). By the proof of Mal'cev-Neumann [cf. \cite{Passman}, Lemma 2.10 p.599-601], $(\widehat T)^+$ is well-ordered and every element belongs to at most finitely many sets $(\widehat T)^n$. Thus we obtain 
$$\lambda^{-1}=1+\sum_{t\in(\widehat T)^+}\alpha_tt\ , \ \alpha_t\in \Q[G_m],$$
thus $\lambda^{-1}\in \Q[G_m]\langle G_m\rangle$. $\qed$

\medskip
The interest of $S(G)$ lies in the following

\begin{proposition} 
\label{goingdown} The projection $\pi_\lambda:{\rm supp}(\lambda)\to G/G_n$ has finite fibers, thus there is a well-defined morphism 
$$\lambda\in \Q[G_n]\langle G/G_n\rangle\mapsto\bar\lambda\in \Q\langle G/G_n\rangle$$
which extends the natural morphism $\Q[G]\to \Q[G/G_n]$.
\end{proposition}

\noindent{\it Proof.}  Writing $\lambda=\displaystyle\sum_{t\in T}\lambda_tt$, we have ${\rm im}(\pi_\lambda)=p(T)$ where $p$ is the projection $T_n\to G/G_n$, which is bijective and increasing since $T_n$ is a transversal. Then $\pi_\lambda^{-1}(\{p(t)\})={\rm supp}(\lambda_t)$, which is finite.

\smallskip To prove "thus", define 
$$\overline\lambda=\sum_{t\in T}\lambda_tp(t)=\sum_{g\in p(T)}\lambda_{p^{-1}(g)}g\in\Q[[G/G_n]].$$ Its support is $\displaystyle\bigcup_{g\in p(T)}{\rm supp}(\lambda_{p^{-1}(g)})$: it is well-ordered since $p(T)$ is ordered as $T$ and the ${\rm supp}(\lambda_t)$ are finite.

 \begin{corollary}
\label{goingup} If $\lambda\in \Z[G_n]\langle G/G_n\rangle$ and $\overline\lambda\in \Z[G/G_n]_{\overline u}$, then $\lambda\in \Z[G]_u$.

\end{corollary}

\noindent{\it Proof.} Let $c\in \R$. By hypothesis, ${\rm supp}(\overline\lambda)\cap\{u<c\}$ is finite. Since ${\rm supp}(\lambda)\to G/G_m$ has finite fibers, ${\rm supp}(\lambda)\cap\{u<c\}$ is finite, thus $\lambda\in \Z[G]_u$.

\section{Proof of Theorem \ref{mixedfinitedetectability}}
\label{proofofmixedfinitedetectability}

By Remark \ref{finiteindex} (3), it suffices to treat the case when $G$ is RTFN. We define the order on $G$, $\Q\langle G\rangle$, $S(G)$ and $\Q[G_m]\langle G/G_m\rangle$ as in the previous section. We assume  ($u(x)>0\Rightarrow x>1$), thus $\Z[G]_u \subset \Q\langle G\rangle$.

\medskip
Let $A  \in {\rm M}_n(\Z[G])$ such that every image in ${\rm M}_n(\Z[G/(H\cap\ker u]_{\overline u})$ (for $H\triangleleft_{f.i.} G$) is invertible. We want to prove that $A$ is invertible in ${\rm M}_n(\Z[G]_u)$. 

\medskip
1) We first prove that $A$ is invertible in ${\rm M}_n(S(G))$. Assume the contrary, then since $S(G) $ is a division ring, there exists 
$L  \in (S(G))^n\setminus\{0\}$ such that $LA = 0$. By Proposition \ref{goingdown}, for $N$ large enough the image 
$\bar L\in (\Q\langle G/G_N\rangle)^n$ is well-defined and nonzero, and we have $\bar L \bar A = 0$, where $\bar A$ is the image of $A$ in ${\rm M}_n(\Q\langle G/G_N\rangle)$. 

\smallskip
Since ${\rm supp}(\bar L)\subset G/G_N$ is finite, we find $k\in\N^*$ such that $m\bar L\in  (\Z\langle G/G_N\rangle)^n$.
Thus $\bar A$ is not invertible in ${\rm M}_n(\Z[G/G_N ])$ and a fortiori in ${\rm M}_n(\Z[G/G_N ]_{\overline u})$. Since $G/G_N$ is nilpotent, Theorem \ref{theorem7.3fornilpotent} (2) implies that for some subgroup $K\triangleleft_{f.i.} G/G_N$,
the image of $\bar A$ is not  invertible in $${\rm M}_n(\Z[(G/G_N)/(K\cap\ker u)]_{\overline u}).$$ 
We have $K=H/H_N$ with $H\triangleleft_{f.i}G$, and there is a natural isomorphism
$$(G/G_N)/((H/H_N)\cap\ker u)\approx G/(H\cap\ker u).$$
Thus the image of $A$ is not invertible in ${\rm M}_n(\Z[(G/(G\cap\ker u)]_{\overline u})$,  contradiction. 

\medskip
2) Let $B$ be the inverse of $A$ in ${\rm M}_n(S(G))$, and let $N_0$ be such that $B\in{\rm M}_n(\Q[G_{N_0}]\langle G/G_{N_0}\rangle)$. By Proposition \ref{goingdown}, for every $N\ge N_0$, $B$ has a well-defined image $\bar B\in {\rm M}_N(\Q\langle G/G_N\rangle )$, and $\bar A  \bar B = {\rm I}_n =\bar B  \bar A$. Thus $\bar B$ is the inverse of $\bar A\in {\rm M}_n(\Q\langle G/G_N \rangle )$. Since $\bar A$ is invertible already in ${\rm M}_n(\Z[G/G_N]_{\overline u})$, we have $\bar B\in {\rm M}_n(\Z[G/G_N]_{\overline u})$. 

\smallskip
Since this is true for all $N\ge N_0$,  $B$ has integer coefficients, ie $B\in{\rm M}_n(\Z[G_{N_0}]\langle G/G_{N_0}\rangle)$. Finally, its image in ${\rm M}_n(\Z\langle G/G_{N_0}\rangle)$ belongs to ${\rm M}_n(\Z[G/G_{N_0}]_{\overline u})$, thus by  
Corollary \ref{goingup}, we have $B  \in {\rm M}_n(\Z[G]_u)$. $\qed$



\begin{thebibliography}{10}

\bibitem[Agol 2014]{Agol 2014} I. Agol, {\it The virtual Haken conjecture}, with an appendix by I. Agol, D. Groves, J. Manning, Doc. Math. 18 (2013), 1045--1087.


\bibitem[Bieri-Neumann-Strebel 1987]{Bieri-Neumann-Strebel 1987} R. Bieri, W.D. Neumann, R. Strebel, {\it A geometric invariant of discrete groups}, Invent. Math. 90 (1987), 451--477.

\bibitem[Chatters 1984]{Chatters 1984}, A.W. Chatters, {\it Noncommutative unique factorization domains}, Math. Proc. Cambridge.Philos. Soc. 95 (1984), 49--54.

\bibitem[Eisenbud]{Eisenbud}  D. Eisenbud, {\it Commutative Algebra. With a View Towards Algebraic Geometry}. Springer Grad. Texts in Math. 150, 1995.

\bibitem[Eizenbud-Lichtman 1987]{Eizenbud-Lichtman 1987} A. Eizenbud, A.I. Lichtman, {\it On embedding of group rings of residually torsion-free nilpotent groups into division rings}, Trans. Amer. Math. Soc. 299 (1987 ), 373--386.

  \bibitem[Friedl-Vidussi 2008]{Friedl-Vidussi 2008} S. Friedl and  S. Vidussi,
 {\it Twisted Alexander polynomials and symplectic structures}, American J. Math. 130 (2008), 455--484.
 

 \bibitem[Friedl-Vidussi 2011]{Friedl-Vidussi 2011}
S. Friedl and  S. Vidussi,
 {\it A survey of twisted Alexander polynomials}, The Mathematics of Knots: Theory and Application (Contributions in Mathematical and Computational Sciences 1, Springer 2011, pp. 45-94.
 
 \bibitem[Friedl-Vidussi 2013]{Friedl-Vidussi 2013}{\it A vanishing theorem for twisted Alexander polynomials with applications to symplectic $4$-manifolds}, J. Eur. Math. Soc. 15 (2013), 2027--2041.

 \bibitem[Goldie 1958]{Goldie 1958} A.W. Goldie, {\it The structure of prime rings under ascending chain conditions}, Proc. London Math. Soc. 8 (1958), 589-608.

 \bibitem[Hall 1954]{Hall 1954} P. Hall, {\it Finiteness conditions for soluble groups}, Proc. London Math. Soc. 4 (1954), 419--436.

 \bibitem[Hall 1959]{Hall 1959} P. Hall, {\it On the finiteness of certain soluble groups}, Proc. London Math. Soc. 9 (1959), 595--622.
 
 \bibitem[Kargapolov-Merzliakov]{Kargapolov-Merzliakov} M.I. Kargapolov, Ju.I. Merzliakov, Fundamentals of the theory of groups. Translated from the second Russian edition by Robert G. Burns. Springer Grad. Texts in Math. 62, 1979.
 
 
\bibitem[Kielak 2020]{Kielak 2020} D. Kielak, {\it Residually finite rationally solvable groups and virtual fibering}, J. Amer. Maths. Soc. 3 (2020), 451--486.

\bibitem[Koberda 2013]{Koberda 2013} T. Koberda, {\it Residual properties of fibered and hyperbolic manifolds}, Topol. Appl. 160 (2013), 857-886.

\bibitem[Kochloukova 2006] {Kochloukova 2006} D.H Kochloukova, {\it Some Novikov rings that are von Neumann finite}, Comment. Math. Helv. 81 (2006), 931--943.

\bibitem[Lam]{Lam} T.Y. Lam, {\it Lectures on Modules and Rings}, Springer Grad. Texts in Math. 189, 1998.


\bibitem[Milnor 1968] {Milnor 1968} J. Milnor, {\it Infinite cyclic coverings}, In: Conference on the Topology of Manifolds, E. Lansing, 1967, p.115-133. Prindle, Weber and Schmidt,1968. Also in the Collected Papers, Volume 2 p.71-95.



\bibitem[Novikov 1981]{Novikov 1981} S.P. Novikov, {\it Multivalued functions and functionals. An analogue of the Morse theory}, Soviet Math. Doklady 24 (1981), 222--226.

\bibitem[Passman]{Passman} D.S. Passman,  {\it The algebraic structure of group rings}, Pure and Applied Mathematics, Wiley, 1977.

\bibitem[Robinson]{Robinson} D.J.S. Robinson, {\it A Course in the Theory of Groups. Second Edition}. Springer Grad. Texts in Math. 80, 1995


\bibitem[Sikorav 1987]{Sikorav 1987} J.-C. Sikorav, {\it  Homologie de Novikov assoc\'ee \`a une classe de cohomologie de degr\'e un}, in: Th\`ese d'État, Universit\'e Paris-Sud (Orsay), 1987.

\bibitem[Sikorav 2017]{Sikorav 2017} J.-C. Sikorav, {\it Novikov homology}, http://perso.ens-lyon.fr/ jean-claude.sikorav/ textes.html.

\bibitem[Stallings 1962]{Stallings 1962} J. Stallings, {\it On fibering certain 3-manifolds}, in {\it Topology of 3-manifolds and related topics}, pp. 95Ð100 Prentice-Hall 1962.

\bibitem[Thurston 1986]{Thurston 1986} W.P. Thurston, {\it A norm for the homology of $3$-manifolds}, in
Mem. Amer. Math. Soc. 59 (1986), no. 339, i--vi and 99--130.


\bibitem[Tischler 1970]{Tischler 1970}
D. Tischler, {\it On fibering certain foliated manifolds over $S\sp{1}$}, Topology 9, 153-154 (1970)



\end{thebibliography}
\end{document}